\documentclass[leqno, 11pt]{amsart}
\usepackage{graphicx,amssymb,epstopdf}
\usepackage[all]{xy}

\newtheorem{thmAB}{Theorem}

\newtheorem{corAB}[thmAB]{Corollary}
\newtheorem{theorem}{Theorem}[section]
\newtheorem{prop}[theorem]{Proposition}
\newtheorem{lemma}[theorem]{Lemma}
\newtheorem{definition}[theorem]{Definition}
\newtheorem{remark}[theorem]{Remark}
\newtheorem{deflemma}[theorem]{Definition and Lemma}

\newcommand{\A}{\mathcal{A}}

\newcommand{\Nn}{\mathcal{N}}
\newcommand{\R}{\mathbb{R}}

\newcommand{\Q}{\mathbb{Q}}
\newcommand{\Z}{\mathbb{Z}}

\newcommand{\image}{\textup{image}}
\newcommand{\length}{\mathcal{L}}
\newcommand{\gra}{\textup{graph}}

\newcommand{\Crit}{\textup{Crit}}
\renewcommand{\span}{\textup{span}}
\newcommand{\id}{\textup{id}}

\newcommand{\st}{ \ : \ }
\newcommand{\Dp}{\Delta_\phi}

\newcommand{\Sp}{S_\phi}
\newcommand{\Hp}{H^\phi}
\newcommand{\Fp}{F^\phi}
\newcommand{\U}{\mathcal{U}}
\newcommand{\V}{\mathcal{V}}
\newcommand{\M}{\mathcal{M}}

\newcommand{\mo}{(M,\omega)}
\newcommand{\om}{\omega}
\newcommand{\Om}{\Omega}
\newcommand{\xm}{x^-}
\newcommand{\xp}{x^+}
\newcommand{\xx}{[x,\widetilde x]}
\newcommand{\yy}{[y,\widetilde y]}
\newcommand{\y}{\widetilde y}
\newcommand{\x}{\widetilde x}
\newcommand{\el}{[l,\widetilde l]}
\newcommand{\xxx}{[x',\widetilde {x'}]}
\newcommand{\zzz}{[z',\widetilde {z'}]}
\newcommand{\yyy}{[y',\widetilde {y'}]}
\newcommand{\zz}{[z,\widetilde z]}

\newcommand{\ls}{[L]\cap^\circ_s(\cdot)}
\newcommand{\lone}{[L]\cap^\U(\cdot)}

\newcommand{\osc}{\textup{osc}}

\newcommand{\dudtau}{\frac{\partial u}{\partial \tau}}
\newcommand{\dudt}{\frac{\partial u}{\partial t}}
\newcommand{\dvdt}{\frac{\partial v}{\partial t}}
\newcommand{\dvdtau}{\frac{\partial v}{\partial \tau}}

\numberwithin{equation}{section}
\numberwithin{figure}{section}

\title[Length Minimizing paths in Ham$(M,\omega)$]
{Length minimizing paths in the Hamiltonian diffeomorphism group}
\author{Peter Spaeth}
\address{Korea Institute for Advanced Study, Seoul, Korea}
\email{spaeth@kias.re.kr}



\date{\today}

\begin{document}
\begin{abstract}
On any closed symplectic manifold we construct a path-connected 
neighborhood of the identity in the Hamiltonian diffeomorphism group with 
the property that each Hamiltonian diffeomorphism in this neighborhood 
admits a Hofer and spectral length minimizing path to the identity.  This 
neighborhood is open in the $C^1-$topology.  The construction utilizes a 
continuation argument and chain level result in the Floer theory of 
Lagrangian intersections.
\end{abstract}
\maketitle

\section{Introduction}

First, we prove a chain level 
result in the Floer theory of Lagrangian intersections, Theorem 
\ref{thm:A}, that relates the local Floer cap product to the so-called 
thin part of the Floer cap product.  The proof relies upon Chekanov's 
version \cite{chekanov:li98} of the Floer continuation argument and a 
generalized definition of chain homotopy, called a $\lambda-$homotopy.  
Together with an algebraic result, Proposition \ref{prop:A}, this provides 
the existence of a certain 
pseudo-holomorphic curve with prescribed asymptotics. See Proposition 
\ref{prop:local_curve}.  The existence of 
this curve yields numerical identities among certain 
non-degenerate Hamiltonian diffeomorphisms.  See Theorem \ref{thm:B}. 
Consequently we achieve the main goal of the article, namely that the 
Hamiltonian group of any closed symplectic manifold exhibits a local 
flatness property under the Hofer and spectral norms. See Corollary C and 
Corollary D.  Under the Hofer norm Bialy and Polterovich \cite{bialy:go94} 
observe this for a neighborhood of the identity of $C^1-$small 
Hamiltonians 
in the group of compactly supported Hamiltonian diffeomorphisms of 
$\R^{2n}$.  Lalonde and McDuff \cite{lalonde:hl95} and McDuff 
\cite{mcduff:gv02} also obtain the $C^1-$small flatness on an arbitrary 
symplectic manifold $\mo$ under the Hofer norm.  Recently Y.-G. Oh 
\cite{oh:si05} proves the $C^1-$small flatness for the Hofer and spectral 
norms.

The neighborhood of the identity we construct is likely larger than has 
appeared in the literature. It consists of all engulfable 
Hamiltonians for which the oscillation of the generating function is 
smaller than half the minimal energy of a non-constant pseudo-holomorphic 
sphere in $M$.  This is an open set in the $C^1-$topology. Engulfability 
is defined in Definition \ref{def:engulfable}.

\subsection{Statement of Theorem \ref{thm:A} }

Let $(P,\omega)$ be any tame symplectic manifold (not necessarily closed).  
This means that there is an $\omega-$tame almost complex structure $J_0$ 
on $P$ such that the corresponding metric $g_{J_0}$ has bounded curvature 
and whose injectivity radius is bounded from zero from below.  Examples 
include all closed symplectic manifolds, the cotangent bundle of any 
closed manifold with the canonical symplectic structure and $\R^{2n}$ with 
its standard symplectic structure.

Let $L$ be a connected, closed Lagrangian submanifold of $P$ and fix a 
Darboux-Weinstein neighborhood $\U$ of $L$. Denote by $\sigma(L;P,\omega)$ 
the minimal area of any non-constant pseudo holomorphic sphere or disk 
with boundary on $L$ in $P.$ Let $H:P\times [0,1] \to \R$ be a smooth 
Hamiltonian function.  Assume 
that $H$ is normalized, which means that for some compact set $K \subset 
P$,
$$
\textup{supp}\ H(\cdot, t) \subset K, \ \textup{for all} \ t \in [0,1].
$$
In addition assume that $H$ is non-degenerate; that is
$$\phi_H^1(L) \pitchfork L.$$
Here $\phi_H^1$ denotes the time-one map of Hamilton's equation
\begin{equation}\label{eq:hamiltons_ode}
\dot z = X_{H}(z,t), \ z(0), z(1) \in L
\end{equation}
where the Hamiltonian vector field $X_{H}$ is defined by the equation
$$
X_{H} \lrcorner \ \omega = dH.
$$

Now suppose $s \in \R$ satisfies $0<s<<1$ and we now consider the 
Hamiltonian function $s\cdot H:P\times [0,1]\to \R.$  
We recall the set-up to the Floer theory $HF(L,sH;\Z_2)$ beginning with 
the action functional.  Fix a reference path $l_0$ on $L$.  The action 
functional $\A_{sH}$ is defined on the relative analogue of Hofer and 
Salamon's \cite{hofer:fh95} Novikov covering 
$\widetilde{\Om}_0(L,l_0;P,\om) = \widetilde {\Om}_0 
\stackrel{\pi}{\to} \Om_0$ over the space $\Omega_0$ of smooth 
contractible paths beginning and ending on $L$
$$
\Omega_0 = \{ l: [0,1] \to P \ | \ l(0), l(1) \in L, \ \textup{$l$ 
homotopic rel $L$ to $l_0$}  \}.
$$
$\widetilde{\Om}_0$ consists of equivalence classes $\el$ where 
the map $\tilde l:[0,1]\times[0,1] \to P$ provides the homotopy from $l_0$ 
to $l$ 
$$\tilde l (0,t) = l_0(t), \ \tilde l(1,t) = l(t), \ \tilde l(s,0) \in L \ 
\textup{and} \ \tilde l(s,1) \in L.
$$
Two pairs 
$(l_1,\widetilde{l_1})$ and $(l_2,\widetilde{l_2})$ are by definition 
equivalent if both $l_1 = l_2$ and the connected sum along $l$, 
$\widetilde{l_1} \# \overline{\widetilde{l_2}}$, 
with the reversed orientation on $\widetilde{l_2}$ satisfies
$$ 
\mu_L(\widetilde{l_1} \# \overline{\widetilde{l_2}}) = 0 = 
I_\omega(\widetilde{l_1} \# \overline{\widetilde{l_2}}),
$$
where $\mu_L$ denotes the Maslov index and $I_\omega$ the symplectic 
action homomorphism on the relative homotopy group $\pi_2(P,L)$.  The 
quotient group
$$
\Gamma = \frac{\pi_2(P,L)}{\ker \mu_L \cap \ker I_\omega}
$$
acts on elements $\el\in \widetilde{\Omega}$ by gluing a disk along 
the reference path $l_0$.  Moreover $\Gamma$ is a free abelian group.

The action functional is defined as follows
$$
\A_{sH}(\el) = -\int_{[0,1]^2} \tilde l^*\omega - \int_0^1 H(l(t),t) \ dt.
$$
Critical points $\zz$ of the action functional consist of pairs
$$
\Crit(\A_{sH}) = \{ \zz \in \widetilde{\Omega}_0 \  | \ z \ 
\textup{satisfies} \ (\ref{eq:hamiltons_ode}) \}.
$$
Typically we will denote elements of the covering space $\el$ or $[l,w]$ 
while critical points will be denoted $\zz$ or $[z,w].$ Denote 
$CF_*(L,sH)$ the vector space of Floer chains generated over $\Z_2$.  The 
solutions to Equation (\ref{eq:hamiltons_ode}) are in one to one 
correspondance with the points of intersection of $L$ and 
$\phi_{sH}^1(L).$

Fix a time-dependent almost complex structure $J=J_t$ on $P$.  The 
relative Floer 
differential is defined by counting finite energy solutions 
$u:\R \times [0,1] \to P$ to the equation
\begin{equation}\label{eq:open_floer}
\left\{ \begin{array}{l}
	\begin{displaystyle}\dudtau + J(u,t) \left( \dudt - X_{sH}(u,t) 
	\right) = 0 \end{displaystyle}\\
	u(\tau,0), u(\tau,1) \in L\\
	\end{array}
\right.
\end{equation}
where the energy is defined to be 
$$
E_J(u):= \int_{-\infty}^{+\infty} \left| \dudtau \right|^2 \ d\tau = 
\int_{-\infty}^{+\infty} \int_0^1 \omega \left( \dudtau, J\dudtau  
\right) \ dt d\tau
$$
with respect to the metric $g_J = \omega(\cdot, J \cdot).$
The energy of a solution $u$ to Equation (\ref{eq:open_floer}) can also be 
written
$$\frac{1}{2} \int_{-\infty}^{+\infty} \left| \dudtau \right| ^2 + \left| 
\dudt - X_H(u,t) \right|^2 \ d\tau.$$
For generic $sH$ the finite energy solutions of Equation 
(\ref{eq:open_floer}) converge exponentially to solutions of Equation 
(\ref{eq:hamiltons_ode}) (see \cite{robbin:ab01}).
We denote the collection of finite energy solutions
$$
\Nn(L,sH;J) = \{ u: \R \times [0,1] \to P \ | \ u \ 
\textup{satisfies} \ (\ref{eq:open_floer}) \ \textup{and} \ E_J(u)<\infty 
\}
$$
and when asymptotic conditions are specified 
$$
\lim_{\tau \to - \infty} u(\tau,\cdot) = \xx,\lim_{\tau \to + \infty} 
u(\tau,\cdot) = \yy \ \textup{with} 
\  \x \# u \simeq \y
$$
we write $\Nn(\xx, \yy ;J)$.  

For a general Lagrangian submanifold $L$ the relative Floer differential 
may be obstructed or not square to zero, however when $s$ is sufficiently 
small and $J$ is sufficiently close to an autonomous almost complex 
structure, the solutions 
to Equation (\ref{eq:open_floer}) are either thin (written $u \in 
\Nn^\circ$) or thick (resp. $u \in \Nn'$) (see also 
Proposition \ref{prop:thinthick} below, Proposition 
4.1 \cite{oh:fc96}, and Lemma 6 \cite{chekanov:li98}).
The thin part of the Floer differential, $\partial^\circ$, is defined by 
counting the thin solutions to Equation (\ref{eq:open_floer}). 
And in fact (again when $s$ is sufficiently small and $J$ is close to an 
autonomous almost complex structure)
$$\partial^\circ \circ \partial^\circ = 0, \ 
H_*(\partial^\circ) = HF^\circ_*(L,sH;\Z_2) \cong 
H_*(L;\Z_2).
$$
See Lemma 7, \cite{chekanov:li98} and Theorem 4.7, \cite{oh:fc96}.

The Floer cap product is defined by counting solutions to Equation 
(\ref{eq:open_floer}) passing through a marked point $p\in L$.
The collection of these maps we denote
$$
\Nn(L,sH;J;p) = \{ u \in \Nn(L,sH;J) \ | \ u(0,0) = p \in L \}.
$$
Again there is a thick and thin decomposition and the thin Floer cap 
product is defined by counting thin elements $u \in 
\Nn^\circ(L,sH;J;p)$.

The thin part of the Floer cap product, for generic $p \in L$, descends to 
an isomorphism on $HF^\circ_*(L,sH;\Z_2)$.  See Lemma 7.4, \cite{oh:st99}.

Next we could also define a local Floer cap product by counting those 
curves $u\in \Nn(L,sH;J;p)$ as above whose image lies in the Darboux 
neighborhood $\U$ of $L$.  These local curves we denote 
$\Nn^\U(L,sH;J,p)$.  While we don't know if this local Floer cap product 
is an isomorphism, it is related to the thin cap action.

Indeed let 
\begin{eqnarray*}
&&\Phi^+: CF_*(L,sH;\Z_2) \to CF_*(L,1\cdot H;\Z_2)\\
&&\Phi^-: CF_*(L,1\cdot H;\Z_2) \to CF_*(L,sH;\Z_2)
\end{eqnarray*}
be the usual continuation maps in Floer theory defined by the Floer 
continuation equation.  (See Section 3 below.) Recall that in 
\cite{chekanov:li98} Chekanov proves that if the Hofer length, 
$||H||$, of H is less than $\sigma(L;P,\omega)$, then for $s$ sufficiently 
small and generic $H$, the identity map and the composition 
$\Phi^-\circ \Phi^+$
$$
\id, \Phi^-\circ \Phi^+ :CF_*(L,sH;\Z_2) \to CF_*(L,sH;\Z_2)
$$
are $\lambda-$homotopic.  See Section \ref{section:algebra} for a review 
of $\lambda-$homotopy. 

This motivated us to prove that the maps $[L]\cap^\circ_s(\cdot)$ and 
$\Phi^-\circ [L]\cap^\U(\cdot) \circ \Phi^+$
$$
[L]\cap^\circ_s(\cdot),\Phi^-\circ \left( [L]\cap^\U(\cdot) \right) \circ 
\Phi^+: CF_s(L,sH;\Z_2) \to CF_s(L,sH;\Z_2)
$$
are $\lambda-$homotopic under a slightly stronger assumption on $||H||$. 

\begin{thmAB}\label{thm:A}
Suppose that the Hamiltonian diffeomorphism $\phi = \phi^1_H$ is 
non-degenerate in the sense that $\phi^1_H(L)$ meets $L$ transversally, 
$\phi^1_H(L) \pitchfork L,$ 
and satisfies
\begin{equation}\label{eq:thmA_assump}
||H|| < \frac{\sigma(L;P,\omega)}{2},
\end{equation}
where $||H||$ is Hofer's length of $H$.  For $s$ sufficiently small the 
thin Floer cap product
$$
[L]\cap_s^\circ(\cdot) : CF_*(L,sH) \to CF_*(L,sH)
$$
is $\lambda-$homotopic to the composition
$$
\Phi^- \circ [L]\cap^\U(\cdot) \circ \Phi^+: CF_*(L,sH) \to CF_*(L,sH).
$$
\end{thmAB}
The proof of Theorem \ref{thm:A} is the subject of Section 3.

\subsection{Theorem \ref{thm:B} and its consequences}
Let $(M,\omega)$ be a closed symplectic manifold of dimension $2n$ and 
$H:M\times [0,1] \to \R$ be any smooth Hamiltonian function.  
The Hamiltonian vector field $X_H$ generates a Hamiltonian flow 
$\phi^t_H$; that is,
$$\frac{d}{dt} \phi^t_H = X_H\circ \phi^t_H, \quad \phi^0_H = \id.$$
We abbreviate this by writing $H\mapsto \phi^t_H$ or $H \mapsto 
\phi$ if we are interested in the time-one map. 
Also we denote the isotopy $t\in [0,1] \mapsto \phi^t_H$ by $\phi_H.$

Given a second Hamiltonian function $G:M\times [0,1] \to \R$ with $G 
\mapsto \psi_G^t$, the 
chain rule $$\frac{d}{dt} \left( \phi^t_H \circ \psi^t_G \right) = 
\dot {\left(\phi^t_H \right)} + T\phi^t_H \dot{\left( \psi^t_H \right)} $$
implies that the Hamiltonian function 
$$H\# G(x,t) := H(x,t) + G((\phi^t_H)^{-1}(x), t)$$ 
generates the composed isotopy $t \mapsto \phi^t_H \circ \psi^t_G$
and that the Hamiltonian function 
$$\overline {H}(x,t):= -H(\phi^t_H x, t)$$  
generates the the inverse flow $t \mapsto (\phi^t_H)^{-1}$.
In other words
\begin{equation}\label{eq:calculus_hamiltonians}
	\left\{
	\begin{array}{l}
	H\# G:= H + G\circ \left(\phi^t_H\right)^{-1} \mapsto \phi^t_H 
	\circ \psi^t_H \\
	\overline{H}=-H\circ \phi^t_H \mapsto \left(\phi^t_H\right)^{-1}.
	\end{array}
	\right. 
\end{equation}

When the symplectic manifold is closed, a Hamiltonian function $H:M\times 
[0,1] \to \R$ is normalized if its mean value is zero for each time $t$ 
$$
\int_M H(x,t) d\mu(x) = 0 \ \textup{for all} \ t \in [0,1],
$$
where $d\mu = \omega^n/n!$.  The Hofer length of a Hamiltonian isotopy 
$\phi_H$, which does not depend on the normalization, is the mean 
value (over time) of the oscillation of the Hamiltonian function $H$
$$
\textup{length}(\phi^t_H) = ||H|| = \int_0^1 \osc_{x\in M}H(x,t) 
\ dt,
$$
where for each $t \in [0,1]$
$$\osc_{x\in M}(H(x,t)) = \max_{x\in M}H(x,t) - \min_{x\in M}H(x,t).$$
Then the Hofer norm of the Hamiltonian diffeomorphism
$\phi$ is
defined to be
$$
||\phi||_{Hofer} = \inf_{H\mapsto \phi} ||H||.
$$
\begin{theorem}[Hofer \cite{hofer:ot90}, Polterovich 
\cite{polterovich:sd93}, Lalonde and McDuff \cite{lalonde:tg95}]The 
function $||\cdot||_{Hofer}: Ham \mo \to \R_+$ satisfies the following 
properties. Let $\phi, \psi \in Ham \mo$.
\begin{enumerate}
\item $||\theta\ \! \phi \ \! \theta^{-1}|| = ||\phi||$ for any $\theta 
\in Symp 
	\mo$ \hfill (symplectic invariance)
\item $||\phi|| = ||\phi^{-1}||$ \hfill (symmetry)
\item $||\phi\ \! \psi|| \leq ||\phi|| + ||\psi||$ \hfill (triangle 
	inequality)
\item $||\phi|| = 0$ if and only if $\phi$ is the identity \hfill 
	(non-degenracy)
\end{enumerate}
\end{theorem}

Hofer \cite{hofer:ot90} proves $d(\phi,\psi) = ||\phi \circ 
\psi^{-1}||_{Hofer}$ defines a bi-invariant non-degenerate distance on the 
group of compactly supported Hamiltonian diffeomorphisms of $\R^{2n}$ with 
its standard symplectic structure using infinite dimensional variational 
methods.  Viterbo \cite{viterbo:st92} later proves this result with 
generating functions. Polterovich \cite{polterovich:sd93} further 
generalizes this result to all tame, rational symplectic manifolds using 
pseudo-holomorphic curves and finally Lalonde and McDuff 
\cite{lalonde:tg95} prove the non-degeneracy property, also via 
pseudo-holomorphic curves, for any symplectic manifold.  Chekanov 
\cite{chekanov:li98} and Oh \cite{oh:gf97} prove the non-degeneracy 
property for tame symplectic manifolds via Floer theoretical techniques.

The association from a Hamiltonian flow to normalized Hamiltonian 
function is injective.  The normalization is important when studying 
the critical values of the associated symplectic action functional.  For 
instance the action spectrum, $\textup{Spec}(H)$, i.e. the set of critical 
values of the action functional $\A_H$ depends only on the homotopy class 
$[\phi,H]$. See \cite{oh:no05}. And in fact if $(M,\omega)$ is 
symplectically aspherical Schwarz \cite{schwarz:ot00} proves 
$\textup{Spec}(H)$ depends only on the time-one map $\phi$.

Now in order to review the spectral invariants (See \cite{oh:cl02} and 
\cite{oh:si05} for a precise treatment) assume temporarily that the 
Hamiltonian function $H$ is 1-periodic in time, $H:M\times \R / \Z \to 
\R.$ Let $\Omega_0(M)$ denote the space of smooth contractible loops 
$l:S^1\to M$ and $\widetilde \Omega_0(M) \to \Omega_0(M)$ denote Hofer and 
Salamon's \cite{hofer:fh95} covering space on which the symplectic action 
functional is defined.  In the closed string setting, $\widetilde 
\Omega_0(M)$ consists of equivalence classes of pairs $(l,w)$, where $w$ 
is a disk bounding $l$.  Let $w \# \overline{w'}$ denote the sphere 
obtained from gluing two disks $w$ and $w'$ (with the opposite 
orientation) along a common boundary.  Two pairs $(l,w)$ and $(l',w')$ are 
said to be $\Gamma-$equivalent if and only if
$$
l = l' \ \textup{and} \ \int_{S^2}(w\#\overline{w'})^* \omega = 0 = 
\int_{S^2}(w \# \overline{w'})^* c_1(\omega).
$$
The automorphism group of the covering $\widetilde{\Omega}_0(M) \to 
\Omega_0(M)$ 
$$ \Gamma = \frac{\pi_2(M)}{\ker \omega \cap \ker c_1(\omega)}$$
is called the Novikov covering group.

We define the action functional $\A_H: \widetilde{\Omega}_0(M) \to \R$
by setting
$$\A_H([l,w]) = -\int_{D^2} w^*\omega - \int_0^1 H(l(t),t) = dt.$$
The assignment $x\in M \mapsto \phi^t_H(x)$ is a one to one correspondance 
between the set of fixed points, Fix$(\phi)$, of $\phi$ and Per$(X_H)$, 
the set of 1-periodic orbits of the Hamiltonian vector field, $X_H$.
A well known calculation shows that the set of cricital 
points $\Crit(\A_H)$ is given by
$$\Crit(\A_H) = \{ \zz \in \widetilde{\Omega}_0(M)\ | \ z \in 
\textup{Per}(X_H) \}.$$
Any constant loop $x \in \textup{Per}(X_H)$ admits the canonical constant 
bounding disk $\widehat{x}$.

The set of Floer-Novikov chains $CF_*(H)$ consists of all formal sums
\begin{equation}\label{eq:floer_novikov_chain}
\alpha_H = \sum_{\zz \in \Crit \A_H} a_{\zz}\cdot \zz, \ a_{\zz} \in \Q
\end{equation}
such that 
$$
\# \{ \zz \ | \ a_{\zz} \neq 0 \ \textup{and} \ \A_H(\zz) \geq 
r \} < \infty$$
for all $r \in \R.$  If the critical point $\zz$ in Equation 
(\ref{eq:floer_novikov_chain}) has a non-zero coefficient $a_{\zz}$ then 
we say that $\zz$ contributes to the chain $\alpha$ or write $\zz \in 
\alpha$. 

Assume that $H$ is non-degenerate in the Floer-theoretical sense and 
denote the Conley-Zehnder index of a critical point $\zz$ by 
$\mu_H(\zz)$. The set of Floer-Novikov chains of index $n$ is denoted 
$CF_n(H).$

The reader is referred to \cite{oh:si05} for a more detailed account of 
the Floer homology theory of the action functional.  Let $J_0$ denote an 
almost complex structure on $M$ and $J = J(\cdot,t)$ be a  
1-periodic family of almost complex structures with $J(\cdot,0) = J_0$.  
We only mention that the Floer differential is defined by counting finite 
energy solutions to the partial differential equation
\begin{equation}
\left\{ 
	\begin{array}{l}
	u : (\tau,t) \in \R \times S^1 \to u(\tau,t) \in M \\
	\begin{displaystyle}\dudtau + J(u,t) \left( \dudt - X_H(u,t) 
		\right) = 0\end{displaystyle}.\\
	\end{array}
\right.
\end{equation}

Any non-zero quantum cohomology class $a \in QH^*(M)$ determines a 
spectral invariant $\rho(H;a) \in \R$ by the following mini-max procedure. 
Fix a Floer-Novikov \textit{cycle} $\alpha$ Poincar\'e dual, written 
$[\alpha]^\flat = a$, (in a precise sense, See the 'Floer Fundamental 
cycle', \cite{oh:co05}) to $a$.  Let its level, $\lambda_H(\alpha)$, be 
the maximum action of any critical point contributing to the cycle.  Then 
$\rho(H;a)$ is defined to be the infimal level among all cycles $\alpha$ 
dual to $a$.  Great care must be made to show this number is finite, does 
not depend on the various choices involved (such as the almost complex 
structure) and satisfies a short list of axioms.  See Theorem I, 
\cite{oh:co05}.

The non-degenerate spectrality axiom (a theorem by the work of Oh 
\cite{oh:fm04} and Usher \cite{usher:sn07}) asserts that when $H$ is 
non-degenerate, $\rho(H;a)$ is a critical value of the action functional 
$\A_H$. The non-degenerate spectrality axiom implies that $\rho(H;a)$ only 
depends on the homotopy class of $\phi_H$ with fixed endpoints, which in 
turn allows the spectral invariant to be extended to any, not necessarily 
periodic, Hamiltonian function $H:M\times [0,1] \to \R$.

Furthermore $\rho(H;a)$ depends continuously on the Hamiltonian $H$ and so 
extends to any smooth function $H$, non-degenerate or not.  The triangle 
inquality asserts that $$\rho(H\#K;a\cdot b) \leq \rho(H;a) + \rho(K;b)$$ 
for any two Hamiltonian functions $H$ and $K$ and where $a \cdot b$ refers 
to the quantum product of $a$ and $b$.  Having said all of this, the 
normalization axiom then implies
$$
0 = \rho(H\# \overline H; 1) \leq \rho(H;1) + \rho(\overline H;1).
$$
The former denotes the spectral length of the isotopy $\phi_H$,
$$
\textup{length}_\gamma(\phi_H) = \gamma(H) = \rho(H;1) + 
\rho(\overline{H};1).
$$

The time-reversal $t\mapsto 1-t$ provides another useful Hamiltonian 
function $\widetilde H (x,t) := - H(x,1-t)$.  The time-one mapping of the 
vector field $X_{\widetilde H}$ is also $(\phi^1_H)^{-1}$ and the two 
isotopies $\phi^t_{\overline H}$ and $\phi^t_{\widetilde H}$ are homotopic 
rel endpoints.  Thus by the homotopy axiom $\rho(\overline H;1) = 
\rho(\widetilde H;1)$ (See Lemma 5.2, \cite{oh:si05}) and so the spectral 
length may be written as
$$\gamma(H) = \rho(H;1) + \rho(\widetilde H;1).$$
The infimum over all isotopies ending at a given Hamiltonian 
$\phi$ defines the spectral norm
$$\gamma(\phi) := \inf_{H \mapsto \phi} \gamma(H) \ .$$
For the reader's conveniece we recall

\begin{theorem}[Schwarz \cite{schwarz:ot00}, Oh \cite{oh:si05}]
The spectral norm
$\gamma: Ham \mo \to \R_+$ satisfies the following 
properties. Let $\phi, \psi \in Ham \mo$.
\begin{enumerate}
\item $\gamma(\theta \ \! \phi \ \! \theta^{-1}) = \gamma(\phi)$ for any 
$\theta 
	\in Symp \mo$ \hfill (symplectic invariance)
\item $\gamma(\phi) = \gamma(\phi^{-1})$ \hfill (symmetry)
\item $\gamma(\phi \ \! \psi) \leq \gamma(\phi) + \gamma(\psi)$ \hfill 
	(triangle inequality)
\item $\gamma(\phi) = 0$ if and only if $\phi$ is the identity \hfill 
	(non-degenracy)
\item $\gamma(\phi) \leq ||\phi||_{Hofer}$
\end{enumerate}
\end{theorem}

The spectral norm $\gamma$ defines a bi-invariant metric on the 
Hamiltonian group which is a Floer-theoretical refinement of the Hofer 
metric. In particular the spectral norm is a lower bound for the Hofer 
norm.  

Roughly speaking the spectral norm of a Hamiltonian diffeomorphism 
$\phi_H$ is the smallest action difference among homologically essential 
critical values of the action functionals $\A_H$ and $\A_{\widetilde{H}}$ 
corresponding to critical points $[z,w]$ and $[\widetilde z,\widetilde w]$ 
with Conley-Zehnder indices $\mu_H([z,w]) = n$ and 
$\mu_{\widetilde{H}}([\widetilde z,\widetilde w]) = n$.  For the moment we 
change the notation briefly and denote critical values of the action 
functionals $\A_H$ and $\A_{\widetilde{H}}$ by $[z,w]$ and $[\widetilde z, 
\widetilde w]$.  Under Poincar\'e duality such Floer cycles correspond 
with the quantum cohomology class $1$.

To reflect the fact that the spectral metric involves only this cohomology 
class, Oh introduces the homological area $A(\phi;1)$ which is, again 
roughly speaking, the minimal energy of a Floer trajectory connecting the 
maximum and minimum of $H$. The time-reversal applied to the critical 
point $[\widetilde z, \widetilde w] \in \Crit \A_{\widetilde H}$ yields a 
critical point, say $[z',w'] \in \Crit \A_H$ and $\mu_H([z',w'] = -n$.  
Fix a point $q \in M$ and Floer cycles $\alpha_H$ and $\beta_{\widetilde 
H}$ of Conley-Zehnder indices $\mu_H([z,w]) = n$, $\mu_{\widetilde 
H}([\widetilde z, \widetilde w]) = n$.

Now consider the one-marked stable Floer trajectories $u$ satisfying the 
following data
$$u = u_1 \# \cdots \# u_N$$
with each $u_j$ a Floer trajectory with finitely many sphere bubbles 
attached satisfying
\begin{equation}\label{eq:data_area}
\left\{
 \begin{array}{l}
	\begin{displaystyle} \dudtau + J(u,t) \left( \dudt - X_H(u,t) 
		\right) = 0 \end{displaystyle}\\
	u_1(-\infty) = [z,w] \in \alpha_H, \ u_N(+\infty) \in [z',w'] \ 
\textup{with} \ [\widetilde z, \widetilde w] \in \beta_{\widetilde H} \\
	u_j(0,0) = q \ \textup{for some} \ j = 1,2, \ldots, N.
 \end{array}
\right.
\end{equation}
The energy $E_J(u)$ of a curve $u:\R \times S^1 \to M$ satisfying Equation 
(\ref{eq:data_area}) is 
$$\int_{-\infty}^\infty \left| \dudtau \right|^2 \ d\tau$$
with respect to the metric on $g_J = \omega(\cdot, J \cdot).$
The energy of such a $u = u_1 \# \cdots \# u_N$ is the sum
$$\sum_k E_J(u_k)$$
plus the symplectic area of any attached sphere bubbles.  Now, the 
homological area is defined by the following mini-max procedure, 
reminiscent of the definition of $\sigma(M,\omega).$

Begin by defining
\begin{eqnarray} 
&&A(\phi, J_0; J; \alpha_H, \beta_{\widetilde H}; q) = \inf \{E_J(u) \ | \ 
u \ \textup{satisfies} \ (\ref{eq:data_area}) \} \nonumber
\end{eqnarray}
and then set
\begin{eqnarray}
&&A(\phi, J_0; J; 1; q) = \inf \{A(\phi, J_0; J; \alpha_H, 
\beta_{\widetilde H}; q) \ | \ [\alpha_H]^\flat = 1 = 
[\beta_{\widetilde H}]^\flat \}.
\end{eqnarray}
Now define
\begin{eqnarray}
A(\phi, J_0; J;1) &=& \sup \{A(\phi, J_0; J; 1; q) \ | \ q \in M \}, 
	\nonumber \\
A(\phi, J_0;1) &=& \sup \{ A(\phi,J_0;J;1) \ | \ J \} \ 
	\textup{and}\nonumber \\
A(\phi; 1) &=& \sup \{ A(\phi, J_0;1) \ | \ J_0 \}.
\end{eqnarray} 

Oh proves in Theorem C \cite{oh:si05} that any non-degenerate Hamiltonian 
diffeomorphism $\phi$, not necessarily engulfable, satisfies
$$
A(\phi;1) \leq \gamma(\phi) \leq ||\phi||_{Hofer}.
$$

The graph $\Delta_\phi$ of a Hamiltonian diffeomorphism $\phi$ is a 
Lagrangian submanifold of the product $(M \times M, -\omega \oplus 
\omega).$  Let $\U$ be a Darboux neighborhood of the the diagonal and 
$\Phi$ be a Darboux chart
$$\Phi :\U \subset (M \times M, -\omega \oplus \omega) \to \V = \Phi(\U) 
\subset (T^*(\Delta), d\Lambda)$$
where $\Phi^*(d\Lambda) = -\omega \oplus \omega$ and $\Lambda$ is the 
canonical one-form on the cotangent bundle $T^*\Delta.$

Recall that a Hamiltonian function $H$ is called quasi-autonomous 
\cite{bialy:go94} if there exist fixed points $x^+,x^- \in M$ with 
$$
H(x^+,t) = \max_{x \in M}H(x,t) \ \textup{and} \ H(x^-,t) = \min_{x \in 
M}H(x,t).
$$
If $\phi_F$ is a Hofer geodesic, then necessarily $F$ must be 
quasi-autonomous.  See \cite{bialy:go94}, \cite{lalonde:hl95} and 
\cite{ustilovsky:cp96}.
 
Motivated by Laudenbach's work on Lagrangian submanifolds 
\cite{laudenbach:es95}, Oh \cite{oh:si05} introduces the following
\begin{deflemma}[Oh \cite{oh:si05}]\label{def:engulfable}A 
Hamiltonian diffeomorphism $\phi$ is \textup{engulfable} if it is 
$C^0-$small in 
the sense that its graph is contained in a Darboux-Weinstein neighborhood 
of the diagonal and its image under the Darboux chart $\Phi$ within the 
cotangent bundle is the graph of an exact one-form, $d\Sp$,
$$
\Phi(\Dp) = \gra \ d\Sp.
$$
When normalized, $\Sp$ is the unique autonomous generating function of $\phi.$ 
Any $C^1-$small diffeomorphism is engulfable.  Let $\Hp$ denote the 
special Hamiltonian function defined by the equation
$$
\Phi(\Delta_{\phi^t_H}) = \gra\ t \ d\Sp,
$$
where $\Delta_{\phi^t_H}$ is the graph of $\phi^t_{\Hp}$ inside 
$(M\times M, -\omega\oplus \omega).$ The Hamiltonian function $\Hp$ is 
quasi-autonomous and so
$$
\int_0^1 \max_{x\in M}\Hp(x,t) - \min_{x\in M}\Hp(x,t)\ dt = 
\osc(\Sp) = \max \Sp - \min \Sp.
$$
\qed
\end{deflemma}
Of course when $\phi$ is also engulfable there is the additional inequality
$$ A(\phi;1) \leq \gamma(\phi) \leq ||\phi||_{Hofer}\leq 
||\Hp|| = \osc(\Sp).$$
 
In the setting of the periodic Floer homology there are thick and thin 
decompositions of the Floer boundary map and pants product (Propositions 
9.1 and 9.3, \cite{oh:si05}) and in the case of $C^1-$small 
diffeomorphisms these thick and thin
decompositions yield the existence of a thin
curve $u$ satisfying Equation (\ref{eq:data_area}) with 
\begin{equation}\label{eq:goal}
\left\{
\begin{array}{l}
	u_1(-\infty) = [\xm, \widetilde{\xm}], u_N(+\infty) = [\xp, 
\widetilde{\xp}]\\
	E_J(u) = \osc(\Sp)
\end{array}
\right.
\end{equation}
which in turn, after some further analysis, implies $A(\phi;1) \geq 
\osc(\Sp)$.  As a result 
$$
A(\phi;1) = \gamma(\phi) = ||\phi||_{Hofer} = \osc(\Sp).
$$ 
See Theorem F and Proposition 9.6, \cite{oh:si05}.

However, when the 
diffeomorphism  $\phi$ is engulfable the thick and thin approach does not 
directly apply and so we must use Theorem \ref{thm:A} to produce the 
existence of a local curve $u$ satisfying (\ref{eq:data_area}) and 
(\ref{eq:goal}).  See Proposition \ref{prop:local_curve}.  This is the 
main step in the proof of
\begin{thmAB}\label{thm:B}
Let $\phi$ be a non-degenerate, engulfable Hamiltonian 
diffeomorphism and assume that the 
oscillation of its generating function $\Sp$ satisfies
\begin{equation}\label{eqn:generating}
\osc(\Sp) < \frac{\sigma (M,\omega)}{2} \cdot
\end{equation}For such a diffeomorphism the following string of equalities 
holds.
\begin{equation}A(\phi;1) =\gamma(\phi) = ||\phi||_{Hofer} = 
||H^\phi|| = \osc(\Sp)
\end{equation}
\end{thmAB}
It is not known whether the homological area depends continuously on 
the Hamiltonian $\phi$ and so Theorem B remains open for degenerate 
$\phi.$  Nevertheless, Theorem B is sufficient to prove 
\begin{corAB}\label{corAB:corA}
Let $\phi$ be any engulfable Hamiltonian, possibly degenerate, 
whose generating function $\Sp$ satisfies $$\osc(\Sp) < 
\frac{\sigma\mo}{2}.$$The path $t \to \phi_{H^\phi}^t$ is Hofer and 
spectral length minimizing.  In fact$$\gamma(\phi) = 
||\phi||_{Hofer} = ||H^\phi|| = \osc(\Sp).$$
\end{corAB}

\begin{corAB}\label{corAB:corB}
The collection $\mathcal H$ of all engulfable Hamiltonians satisfying 
Condition (\ref{eqn:generating})
$$
\mathcal H \subset \textup{Ham}\mo
$$
is an open, path-connected neighborhood of the identity in the 
$C^1-$topology such that any element of $\mathcal H$ admits a Hofer 
and spectral length minimizing path to the identity.  The 
neighborhood $\mathcal H$ with either the Hofer or spectral norm is 
isometric to the normed vector space $C_m^\infty(M)$ of mean zero 
functions $f:M\to \R$ with the norm given by the oscillation, $\osc(f) 
= \max f - \min f$.
\end{corAB}

For $C^1-$small Hamiltonians Lalonde and McDuff \cite{lalonde:hl95} and 
McDuff \cite{mcduff:gv02} obtain this result under the Hofer metric, and 
Oh \cite{oh:si05} for both the Hofer and spectral distances.

Suppose the symplectic manifold is symplectically asperical.  In this case 
hypothesis (\ref{eqn:generating}) is not restrictive since there are no 
non-constant pseudo-holomorphic spheres and $\sigma\mo = +\infty$.  It is 
also worth comparing Corollary C to Ostrover's result 
Theorem 1.2, \cite{ostrover:ac03} 
that when $\mo$ is symplectically aspherical, there exist Hamiltonian 
diffeomorphisms (far from the identity) admitting no Hofer length 
minimizing path to the identity.  An example of this kind is also given on 
$S^2$ by Lalonde and McDuff \cite{lalonde:hl95}.

An explicit description of the size of the neighborhood, $\Nn \subset Ham 
\mo$, McDuff constructs is explained in Remark 3.5, \cite{mcduff:gv02}.  
It appears that this description requires some control on both the Hofer 
length and $C^2-$ norm of the generating Hamiltonian function, which we 
avoid.  

In some cases the Darboux neighborhood of the diagonal $\U \subset M 
\times M$ can be quite large.  One may ask if there is a relationship 
between the diameter of the Hamiltonian diffeomorphism group and the size 
of the Darboux neighborhood of the diagonal.

\smallskip 
\begin{center} 
\textsc{Acknowldgements} 
\end{center} 
\smallskip 

This work formed a portion of my Ph. D. thesis.  I wish to thank my 
advisor Yong-Geun Oh for his guidance and support.  Theorem B is motivated 
by a version of that result announced in the preprint \cite{oh:si04}. I 
also wish to thank Joel Robbin and Augustin Banyaga for many encouraging 
discussions.  A final thanks goes to the referee whose comments 
significantly improved the article.

 
\section{$\lambda$-homotopy}\label{section:algebra}

Chekanov \cite{chekanov:li98} defines an algebraic framework to 
prove an existence theorem of Lagrangian intersections.  
We recall the definitions and prove a new algebraic result,
Proposition \ref{prop:A}, needed in the proof of Theorem \ref{thm:B}.

Let $\Gamma$ be a free abelian group and $\lambda: \Gamma \to \R$ be an 
injective homomorphism.  Such a mapping is called a weight function.  
Denote
\begin{equation}
\Gamma^+ = \{ g \in \Gamma | \lambda(g) > 0 \} \qquad
\Gamma^- = \{ g \in \Gamma | \lambda(g) < 0 \}.
\end{equation}

Let $k$ be a commutative ring with nonzero unity and $K=k(\Gamma)$ be the 
group ring over $k$.   The group ring $K$ includes the group $\Gamma$ as a 
subgroup.

Now let $M$ be a $k$-module.  The group $\Gamma$ need not act on $M$, but 
$\Gamma$ does act on $M\otimes K$ via the second factor $g(v\otimes \overline g):= 
v \otimes g(\overline g).$ Furthermore $M\otimes K$ admits the decomposition
\begin{equation}\label{eq:alg_decomp}
M\otimes K = M^+ \oplus M^0 \oplus M^-
\end{equation}
where $M^+ = \Gamma^+(M\otimes K), \quad M^0 = M \ \textup{and} \ M^- =  
\Gamma^-(M\otimes K).$

The above splitting supports two projections
$$p^+:M\otimes K \to M^+ \oplus M^0 \ \textup{and} \ p^-: M\otimes K \to 
	M^0\oplus M^-,
$$
$p^+(v\otimes \overline g) = v \otimes \sum p^+(\alpha_i)g_i$
where 
$$p^+(\alpha_i) =
	\left\{ \begin{array}{ll}
	0 & \textup{if} \ \lambda(g_i) < 0 \\
	\alpha_i & \textup{if}\  \lambda(g_i) \geq 0
		\end{array}
	\right.
$$
and similarly for $p^-.$

Now assume in addition that $(M,\partial)$ is a differential
$k-$module with $\partial \circ \partial = 0.$  Notice that $\partial$ 
extends to $M\otimes K$ by setting $\partial(v \otimes \overline g) = 
(\partial v) \otimes \overline g.$
Also the differential $\partial$ and the projections 
$p^+,p^-$ are commuting operators because $\partial$ acts on $M$ while the 
projections act only on $K$ in $M\otimes K.$
\begin{definition}[Chekanov, \cite{chekanov:li98}] Two $K-$linear maps 
$\psi_1,\psi_2: M\otimes 
K \to M \otimes K$ are called \textup{$\lambda-$homotopic} if there exists 
a $K-$linear map $h: M \otimes K \to M \otimes K$ such that
\begin{equation}\label{eqn:lambda_htpy}
p^+ \left( \psi_1 - \psi_2 + h \circ \partial + \partial \circ h \right) 
p^- = 0.
\end{equation}
\end{definition}
Notice that if $\psi_1,\psi_2$ were chain homotopic maps in the usual sense, 
then Equation (\ref{eqn:lambda_htpy}) holds.  The definition implies that 
$\psi_1$ and $\psi_2$ are chain maps from negative to positive weight. 
The essential algebraic lemma, observed first by Chekanov, from which this 
section really follows is 
\begin{lemma}[Chekanov, \cite{chekanov:li98}]\label{lemma:key_lemma}Let $V$ 
be the homology of $(M,\partial).$  
Set $W=V\otimes K.$  For any nonzero $r \in W$ there exists $g \in \Gamma$ 
for which $g\cdot r = q^0 + q^-$, with $q^0 \in V \subset M^0$ nonzero and 
$q^- \in M^-.$ 
\end{lemma}

\proof We have $r=v\otimes \overline g \in W = V \otimes K,$ with 
$\overline g = 
\sum \alpha_i g_i, \alpha_i \in k, g_i \in \Gamma$.  Let $g=g_j \in 
\Gamma$ be such that 
$$g = g_{max}^{-1} \ \textup{with} \ \lambda(g_{max}) = \max_{\alpha_i 
\neq 0} \lambda(g_i).
$$
The quantity $\max_{\alpha_i \neq 0} \lambda(g_i)$ is called the 
\textit{valuation} of $g$.  Then by the action of $\Gamma$ on $K$ $$g (r) 
= v\otimes \sum_{g_ng_j=g_i}\alpha_n g_i.$$ Finally if 
$\alpha_n \neq 0$ then $\lambda(g_i) = \lambda(g_n) + \lambda(g_j) \leq 
-\lambda(g_j) + \lambda(g_j) =0$ and hence $\lambda(g_i) \leq 0.$ \hfill 
$\qed$

What follows below is the main algebraic result of the current paper.  
An isomorphism on homology can not be $\lambda-$homotopic to the zero 
map provided the isomorphism commutes with the projection $p^+$.

\begin{prop}\label{prop:A}Suppose $(M,\partial)$ is a differential 
$k-$module with non-trivial homology.  Suppose the linear maps
$\psi_1,\psi_2 : M \otimes K \to M\otimes K$
are $\lambda-$homotopic.  If $\psi_1$ descends to an isomorphism 
$\psi_1:H(M,\partial) \stackrel{\cong}{\to} H(M,\partial)$ 
and the induced map commutes with $p^+$ on homology (i.e. up to a boundary 
on the chain level) then $\psi_2 \neq 0.$
\end{prop}

\proof Choose a submodule of $\partial-$cycles $V \subset M = M^0$ 
representing the 
homology $H(M,\partial)$.  Let $W = V \otimes K$.  Note that $W$ is not trivial.  
By Lemma \ref{lemma:key_lemma} for any nonzero $r \in W$ there exists $g 
\in \Gamma$ for which $g\cdot r = q^0 + q^-$, with $q^0 \in V \subset M^0$ 
nonzero and $q^- \in M^-.$

To prove the proposition suppose to the contrary that $\psi_2 = 0$.  Applying equation 
(\ref{eqn:lambda_htpy}) to $q^0 + q^-$ we see
\begin{eqnarray*}
0&=&p^+ \left( \psi_1 - \psi_2 + h \circ \partial + \partial \circ h 
	\right) p^- (q^0+q^-) \\
&=& p^+\left( \psi_1 - \psi_2 + h \circ \partial + \partial \circ h 
	\right)(q^0+q^-) \\
&=&p^+(\psi_1(q^0+q^-)) + p^+(h \circ \partial + \partial \circ h 
	)(q^0+q^-) \\
&=& p^+(\psi_1(q^0+q^-)) + \partial (p^+ h (q^0+q^-))
\end{eqnarray*}

Here we have used the linearity of all the maps, that $q^0 +q^-$ is a 
cycle and that the projection $p^+$ commutes with 
the boundary.  By assumption $p^+$ and $\psi_1$ commute up to a 
boundary, say 
$\partial c$ and so continuing
\begin{eqnarray*}
0 &=& \psi_1(q^0) +\partial c +  \partial (p^+ h (q^0 + q^-))\\
&=& \psi_1(q^0) + \partial \left( c + p^+h(q^0 + q^-) \right)
\end{eqnarray*}
which implies that $\psi_1(q^0)$ is a boundary.  Hence, $\psi_1$ can not 
be an 
isomorphism, a contradiction. \hfill $\qed$

\begin{remark}
\textup{The downward Novikov ring $\Lambda^\downarrow$ is defined to be
$$
\Lambda^\downarrow = \left\{ \sum_i \alpha_i g_i \ | \ (\forall) t \in \R
\ \#\{g_i\ | \ \alpha_i \neq 0 \ \text{and} \ \lambda(g_i) > t\} < \infty \right\}.
$$
Because the valuation of any element in the
downward Novikov ring is finite, Lemma \ref{lemma:key_lemma} continues
to hold with the group ring $K$ replaced with $\Lambda^\downarrow$.  In
other words, the notion of $\lambda-$homotopy extends to this coefficient 
ring.} 
\end{remark} 


\section{The Proof of Theorem \ref{thm:A}} 

Let $L$ be a closed, connected $n-$dimensional Lagrangian submanifold of 
the tame symplectic $(P,\omega).$ Let $\U$ be a Darboux-Weinstein 
neighborhood of $L$.  Fix a reference path $l_0$ 
on $L$.  Given a (normalized) time-dependent Hamiltonian function 
$H:P\times [0,1] \to \R$ the corresponding action functional $\A_H$ is 
real valued on the Novikov covering $\widetilde{\Om}_0(L,l_0;P,\om) = 
\widetilde {\Om}_0 \stackrel{\pi}{\to} \Om_0$ of the space $\Omega_0$ of 
contractible paths beginning and ending on $L$.  $\widetilde{\Om}_0$ 
consists of equivalence classes $\ell$ where the half disk $w$ provides 
the homotopy to $l_0$ and two pairs $(l_1,\widetilde{l_1})$ and 
$(l_2,\widetilde{l_2})$ are equivalent if $l_1 = l_2$ and the Maslov index 
$\mu_L$ and symplectic action $I_\omega$ homomorphisms vanish on the glued 
disk $\widetilde{l_1}\# \overline {\widetilde{l_2}}$.  The Maslov index 
conventions and the Novikov covering are discussed in greater detail in 
\cite{fukaya:li06}, Chapter 2, Sections 2 and 3.

\setlength{\unitlength}{1.25cm}
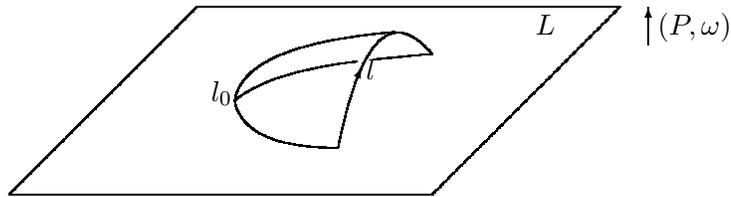
\begin{figure}[!h]
\begin{center}
\begin{picture}(7,5)(-3.5,-2.5)
\qbezier(-3,-1)(-2,0)(-1,1)
\qbezier(1.5,-1)(2.5,0)(3.5,1)
\qbezier(-1,1)(1.25,1)(3.5,1)
\qbezier(-3,-1)(-.75,-1)(1.5,-1)
\qbezier(0.5,-0.5)(.85,1.25)(1.5,0.5)
\qbezier(-0.6,0)(-0.5,-0.5)(0.5,-0.5)
\qbezier(-0.6,0)(-.5,.6)(1.1,0.731)
\qbezier(-0.6,0)(-0.25,.3)(0.7,0.42)
\qbezier(0.875,0.44)(1.1,0.46)(1.5,0.5)
\put(-0.85,0){$l_0$}
\put(2.6,.7){$L$}
\put(3.9,.7){$(P,\omega)$}
\put(3.8,.6){\vector(0,1){0.4}}
\put(0.8,0.225){$l$}
\put(0.7565,0.35){\vector(1,2){0}}
\end{picture}
\caption{Half-disk $\widetilde{l}$ with boundary on $L$}\label{fig:disk}
\end{center}
\end{figure}
The choice of sign is
$$
\A_H(\ell) := -\int \widetilde{l}^*\omega - \int_0^1 H(l(t),t)\ dt.
$$ 
Denote $\Gamma = \Gamma(L,l_0;P,\om)$ the automorphism group of this
cover, which is a free abelian group. 
$$
\Gamma = \frac{\pi_2(P,L)}{\ker (\mu_L) \cap \ker( I_\omega)}
$$
$\Gamma$ acts by gluing a disk: $g \cdot
\ell = [l,\widetilde{l}\# g].$ The symplectic action provides the weight 
homomorphism
$$
\lambda(g):= -\int g^*\omega.
$$

Let $\sigma(L;P,\om)$ be the minimal area of a pseudo-holomorphic sphere
or disk with boundary on $L$ in $P$ and assume for the remainder of this
section that the Hofer length of $H$ satisfies $||H||< \sigma/2.$ Let
$0<s<<1$ and fix the constant $\delta>0$ to satisfy $||H|| < \delta/2 < 
\sigma/2.$

For the parametrized Hamiltonian $sH$ denote the critical points of the
action functional $\A_{sH}$ 
$$
\Crit(\A_{sH}) = CF(L,sH) = CF(s)$$ 
and the vector space of Floer chains
$$
CF_*(s) = CF_*(L,sH;\Z_2) = \span_{\Z_2}(CF(s)).
$$
We will make a number of assumptions on the size of the parameter $s$ in
the course of the proof of Theorem \ref{thm:A}; however, this
parameter does not appear in the statement of Theorem
\ref{thm:B}.  The first of these conditions is that 
\begin{equation} \phi^t_{sH}(L) \subset \U \quad (\forall) \ t\in [0,1].
\end{equation}

Two elements of $CF(s)$ are said to be equivalent (See section 4, 
\cite{chekanov:li98}) if they belong
to the same connected component of the set
$$
\pi^{-1}\{ l\in \Om_0 \st l([0,1])\subset \U \} \subset \widetilde{\Om}_0.
$$

The first assumption on the parameter $s$ guarantees that the set of
equivalence classes is not empty. Fixing any equivalence class
$CF^\circ(s)$ yields the decomposition $CF_*(s) = CF^\circ_*(s) \otimes
K$, where $CF^\circ_*(s) = \span_{\Z_2}CF^\circ(s)$ and $K = \Z_2(\Gamma)$
denotes the group ring.  Observe that for two equivalent critical points 
$\xx, \zz \in CF^\circ(s)$ there is a topological (half) cylinder $C$ with 
boundary 
on $L$ contained in $\U$ such that $\widetilde{x} \# C \simeq 
\widetilde z.$  Hence $\widetilde x$ 
and $\widetilde z$ have equal symplectic area.  Thus for equivalent 
critical points
$$
\A_{sH}(\xx) - \A_{sH}(\zz) = s \left( \int H(z(t),t) - H(x(t), t)\ dt 
\right)
$$
and so if we assume further that $s< \delta/(2||H||+1)$ we will have
\begin{equation}
\xx,\zz \in CF^\circ(s) \implies \A_{sH}(\xx) - \A_{sH}(\zz) < \delta / 2
\end{equation}
The sets of positive and negative weight critical points are denoted
$$
CF^+(L,sH) = \Gamma^+(CF^\circ) \ \textup{and} \ CF^-(L,sH) 
	= \Gamma^-(CF^\circ).
$$

Let $J$ be a (possibly time-dependent) almost complex structure 
compatible to $\omega.$ Again The Floer differential is defined by 
counting finite energy solutions of Equation (\ref{eq:open_floer}).

The full Floer moduli space of finite energy solutions 
$\Nn(L,sH;J,\omega) =\Nn(s)$
ranges over all critical points of $\A_{sH}$.  The finite energy
solutions converge exponentially as $\tau \to \pm \infty$ to paths
$x^\pm:[0,1] \to P$
satisfying Equation (\ref{eq:hamiltons_ode}).

For generic $s$, $\Nn(\xx,\zz;s)$ is a smooth manifold of dimension 
$\mu_L(\xx) - \mu_L(\zz)$, when it is not empty. See \cite{gromov:pc85}.

The length of a Floer trajectory $u \in \Nn(\xx,\zz,s)$ is defined to be 
the action difference at its ends
\begin{equation}
\length(u):= \A_{sH}(\xx) - \A_{sH}(\zz).
\end{equation}

We say that $u$ is \textit{somewhat thin} if $\length(u) \leq ||H|| + 
\delta$, \textit{thin} if $\length(u) \leq \delta$ and \textit{local} if 
$\image(u) \subset \U.$ Curves which exit the Darboux-Weinstein 
neighborhood are called \textit{thick}.  A consequence of Gromov's 
compactness theorem is Oh's thick and thin decomposition theorem that 
implies when $s$ is sufficiently small and $J$ is $C^0-$close to a time 
independent compatible almost complex structure that every somewhat thin 
trajectory is thin and every (somewhat) thin trajectory is local.  (See 
Lemma 6, \cite{chekanov:li98} and Proposition 4.1, \cite{oh:fc96}.) Our 
third assumption on the size of the parameter $s$ is that it is small 
enough to ensure the thick and thin decomposition. We denote the 
collection of thin trajectories $\Nn^\circ(L,sH)$ and local trajectories 
$\Nn^\U(L,sH).$ After dividing by the time-shift, the thin Floer 
differential $\partial^\circ$ is defined by $\Nn^\circ/\R:= 
\widehat{\Nn}^\circ(L,sH)$
$$
\partial^\circ(\xx) = \sum_{\mu_L(\xx) - \mu_L(\zz) = 1} \#_{\Z_2} 
	\widehat{\Nn^\circ}(\xx,\zz,s)\zz
$$
When $s$ is sufficiently small and the almost complex structure is 
close to an autonomous almost complex structure, the 
thin differential satisfies $\partial^\circ \circ \partial ^\circ = 0$ 
and the resulting homology is isomorphic to the singular homology 
of $L$
\begin{equation}
H(CF_*(s), \partial^\circ) := HF^\circ_*(L,sH;\Z_2) \cong H_*(L;\Z_2).
\end{equation}
See Lemma 7, \cite{chekanov:li98} and Theorem 4.7, \cite{oh:fc96}.

The thin Floer cap product $\ls:CF^\circ_*(s) \to CF^\circ_*(s)$ is
defined by counting elements $u \in\Nn^\circ(L,sH)$ passing through a
generic point $p \in L$.
$$\Nn^\circ(L,sH;p) = \{ u \in \Nn^\circ \ | \  u(0,0) = p \}$$ 
$$[L]\cap^\circ_s(\xx) = \sum_{\mu_L(\xx) - \mu_L(\zz)=n} \#_{\Z_2}
\Nn^\circ(\xx,\zz,s;p)\zz$$

The collection of such maps is denoted $\Nn^\circ(L,sH;p)$ (or
$\Nn^\circ(\xx,\zz,s;p)$ or $\Nn^\U(L,sH;p)$, etc.).  Because of the pointed
condition $u(0,0) = p$, the dimension is reduced by $n$, the dimension of
$L$.  Again the thick and thin decomposition holds.

\begin{prop}[Proposition 4.1, \cite{oh:fc96} and 
Lemma 6, \cite{chekanov:li98}]\label{prop:thinthick}
Fix a time-independent compatible almost complex structure $J_0$ on $P$.  
There exists $s'>0$ with the property that for any $0<s<s'$ if 
$||J-J_0||_{C^0(P\times [0,1])}<s$ then any somewhat thin cap trajectory 
is thin and any (somewhat) thin cap trajectory is local.
\end{prop}

Proposition \ref{prop:thinthick} implies that the moduli space $\Nn^\circ(L,sH;p)$ is
isolated, say by the length defined above, within $\Nn(L,sH;p)$.  Passing
to the cotangent bundle, the work
of \cite{floer:mt88}, \cite{oh:st99} and \cite{fukaya:zl97} implies that $\ls$
descends to an isomorphism on $HF^\circ_*(L,sH;\Z_2)$.  (In particular A.
Floer originally proved the injectivity of the cap action.)  The next
assumption on $s$ is that is it small enough to ensure
\begin{equation}\label{eqn:cap_iso}
\ls : HF^\circ_*(L,sH;\Z_2) \stackrel{\cong}{\to} HF^\circ_{*-n}(L,sH;\Z_2).
\end{equation}
The thin cap product shifts the grading by $n$, the dimension of the Lagrangian $L$. 

Setting the parameter $s=1$, the local cap product $\lone$ is defined by the 
local curves in the moduli space $\Nn^\U(L,H;p) := \{ v\in \Nn(L,1\cdot H;p) 
\ | \ 
\textup{Image}(u) \subset \U\}.$
$$[L]\cap^\U(\yy) = \sum_{\mu_L(\yy) - \mu_L(\yyy) = n} \#_{\Z_2} 
	\Nn^\U(\yy,\yyy;p)\yyy$$
Within the Darboux neighborhood of $L$ bubbling-off can not occur and so 
$[L]\cap^\U(\cdot)$ is defined; however, this map need not agree with the 
singular cap product.  We have thus far fallen short in proving that the 
moduli space of local curves is isolated among all curves.  This fact is 
what motivated the use of $\lambda-$homotopy.  For the same reason, when 
$s=1$ the local Floer homology might not be defined.  In particular 
$\partial^\U$ need not square to zero.

The continuation maps $\Phi^+$ and $\Phi^-$ are constructed by introducing 
monotone continuation functions $\rho^+$ and $\rho^-$ of $\tau$ and a 
non-autonomous Floer equation.  We sketch definition of the map $\Phi^+$, 
referring the reader to \cite{chekanov:li98} for the full details.  Let 
$J^{\rho^+}$ be a $\tau$ (and $t$) dependent almost complex structure 
such that is independent of $\tau$, when $|\tau| >> 0.$

Let $\begin{displaystyle}b_+(H) = \int_0^1 \max_{x}H(x,t) \ 
dt\end{displaystyle}$ and $\begin{displaystyle}b_-(H) = -\int_0^1 
\min_xH(x,t) \ dt\end{displaystyle}$.  Let $\rho^+$ be a smooth increasing 
function of $\tau$ so that $\rho^+(\tau) = 1$ for $\tau >> 1$ and 
$\rho^+(\tau) = s$ for $\tau <<-1.$ Consider the continuation equation
\begin{equation}\label{eqn:continuation}
\left\{
\begin{array}{l}
\begin{displaystyle}\dudtau + J^{\rho^+}(u,t) \left( \dudt - X_{\rho^+H}(u,t) \right) = 0\end{displaystyle} \\
u:\R \times [0,1] \to P \\
u(\tau,0), u(\tau,1) \subset L \\
\end{array}
\right.
\end{equation}
and define $$\Nn_{\rho^+} = \{ u:\R \times [0,1] \to P \ |
 \ E_J(u) < + \infty \ \textup{and} \ u \ \textup{satisfies} \ (\ref{eqn:continuation}) \}$$

The length of a continuation trajectory $u^+ \in \Nn_{\rho^+}$ can be 
defined as before (and could be negative) $u^+ \in \Nn_{\rho^+}$ 
$$
\length(u^+) = \A_{sH}(\xx) - \A_{H}(\yy)
$$
where $u^+(-\infty) = \xx\in CF(s)$ and $u^+(+\infty) = \yy\in CF(1)$.  
The continuation trajectories $u^+,u^-$ are called thin if 
\begin{equation}\label{eq:thin_continuation_def}
	\left\{
	\begin{array}{l}
	\length(u^+) \leq (1-s)\cdot b_+(H) + \delta\\
	\length(u^-) \leq (1-s)\cdot b_-(H) + \delta
	\end{array}
	\right. 
\end{equation}
and the collection of all thin continuation trajectories we denote 
$\Nn^\circ_{\rho^+}$ and $\Nn^\circ_{\rho^-}$.Observe that to be thin, a 
continuation trajectory satisfies a weaker inequality than is required for 
a Floer trajectory to be thin.

\begin{lemma}[Lemma 9, \cite{chekanov:li98}]\label{lemma:cont_thin} If 
$u^+ \in \Nn_{\rho^+}$, then $\length(u^+) \geq -(1-s)\cdot b_-(H)$.  If 
$u^- \in \Nn_{\rho^-}$, then $\length(u^-) \geq -(1-s)\cdot b_+(H).$ For 
the thin continuation trajectories $E(u^\pm) \leq ||H|| + \delta$. 
\end{lemma}

Consequently the continuation map $\Phi^+: CF(s) \to CF(1)$ defined by
$$
\Phi^+(\xx) = \sum_{\mu_L(\xx) - \mu_L(\yy)=0} \#_{\Z_2} 
\Nn^\circ_{\rho^+}(\xx,\yy)\yy
$$
is finite. Again we refer to the proof of Lemma 9, \cite{chekanov:li98} 
for 
the details of the finiteness argument.  Notice that a thick and 
thin decomposition is not needed to make 
this definition.

The important point is that these linear maps provide a bridge from the Floer 
complexes corresponding to $0<s<<1$ and $s=1,$
\begin{equation}
\Phi^+: CF_*(L,sH) \to CF_*(L,H) \ \textup{and} \ \Phi^-:CF_*(L,H) \to 
CF_*(L,sH).
\end{equation}
The main result in \cite{chekanov:li98} is that when $s$ is sufficiently 
small and $||H||< \sigma$, there exists a $\lambda-$homotopy $h$ from the 
identity to the composition $\Phi^- \circ \Phi^+.$ As a consequence the 
displacement energy of any Lagrangian submanifold is larger than $\sigma.$ 
Our $\lambda-$homotopy is defined similarly, where we must also take into 
consideration the marked point $p\in L.$ We now provide the details.

Let $R\in [s,+\infty)$ and $\mu_R: \R \to [s, 1]$ be a smooth function 
with all the following properties.  Firstly $0 \leq \mu_R'(\tau)$ for 
$\tau \leq 0$ and $\mu_R'(\tau) \leq 0$ for $\tau \geq 0$ and for $R\geq 
1$
\begin{equation}
 \mu_R(\tau) = \left\{ \begin{array}{ll}
 		1 & -R \leq \tau \leq R\\
		s & |\tau| \geq R+1
	   \end{array}
  \right.
\end{equation}
while when $s\leq R\leq 1,$ $\mu_R = R \cdot \mu_1$.  Let $J^{\mu_R}$ be 
a $\tau-$dependent almost complex structure that is $\tau-$independent 
for $\tau >> 0.$   The $\lambda-$homotopy is 
defined via the finite energy solutions of the non autonomous Floer equation
\begin{equation}
 \left\{
 	\begin{array}{l}
	\begin{displaystyle} \dudtau + J^{\mu_R}(u,t) \left( \dudt - 
		X_{\mu_R(\tau)H}(u,t)\right) = 0\end{displaystyle}\\
	u(\tau,0), u(\tau,1) \subset L\\
	u(0,0) = p \in L
	\end{array}
  \right.
\end{equation}
The collection of all finite energy solutions is denoted $\Nn_{\mu_R}= 
\Nn_{\mu_R}(L,sH;p)$.  The parametrized moduli space is defined to be 
$\Nn_\mu = \{(R,v) \st v \in \Nn_{\mu_R}(L,sH;p) \}$ and by the 
parameterized 
index theorem (See \cite{oh:gf97}) the local dimension of $\Nn_\mu$ is one 
greater than the difference of Maslov indices of the asymptotic limits and 
so defines a map of degree $+1.$ The formula for the pointed 
$\lambda-$homotopy $h^p:CF(s) \to CF(s)$ on generators is
\begin{equation}\label{eqn:homotopy}
h^p(\xx) = \sum_{\yy \in CF^\circ \cup CF^+}\#_{\Z_2} \left\{\Nn_\mu(\xx,\yy;p)\right\}\yy
\end{equation}
where $\xx \in CF^-\cup CF^\circ.$ Because the energy is controlled by the Hofer norm, 
bubbling off does not occur and so the sum is finite.  The induced map is linear over $K = \Z_2(\Gamma)$.

We now verify the $\lambda-$homotopy equation.
\begin{equation}\label{eqn:to_verify}
p^+([L]\cap^\circ_s(\cdot) - \Phi^- \circ [L] \cap^\U(\cdot) \circ \Phi^+ + h^p \circ \partial^\circ + \partial^\circ \circ h^p )p^- = 0
\end{equation}
for all $\xx \in CF^\circ(s).$  We write the left hand side of 
Equation (\ref{eqn:to_verify}) as $$\sum \#_{\Z_2}\mathcal{S}(\xx,\zz)\zz$$
where $\mathcal{S}(\xx,\zz)$ is the union of the following four sets.   
For the remainder of this section each of 
$\xx \in CF^\circ(s)\cup CF^-(s), \zz \in CF^+(s)\cup CF^\circ(s)$ 
and $\yy,\yyy \in CF(1)$ represent generic critical points of the respective action functionals.
\begin{enumerate}
\item A thin cap trajectory $u\in \Nn^\circ(\xx,\zz,s;p)$
\item A triple $(u^-,v,u^+)$ with
$u^+ \in \Nn^\circ_{\rho^+}(\xx,\yy)$, $v \in \Nn^\U(\yy,\yyy;p)$ and $u^- \in \Nn^\circ_{\rho^-}(\yyy,\zz)$
\item A thin unparametized gradient trajectory 
$u \in \widehat{\Nn^\circ}(\xx,\zz,s)$ and element $(R,v)$ with $v \in \Nn_{\mu_R}(\xx,\zz;p)$
\item An element $(R,v)$ with $v \in \Nn_{\mu_R}(\xx,\zz;p)$
and a thin unparametized gradient trajectory $u \in \widehat{\Nn^\circ}(\xx,\zz,s)$
\end{enumerate}

By the Floer gluing theorem each term above lies at the end of the 
one-dimensional portion $\mathcal{Q}$ of $\mathcal{S}$.  In fact 
$\mathcal{Q}$ is compact and the correspondence is one-one.  The 
compactness of $\mathcal{Q}$ follows from the uniform energy estimate 
which is a consequence of the assumption $||H||_{Hofer} < \sigma / 2$ (See 
Lemma 10, \cite{chekanov:li98} and Lemma 2.2, \cite{oh:gf97}.).

On the other hand by the Gromov-Floer compactness theorem (see Lemma 6.2 
\cite{salamon:mt92}) every end of $\mathcal{Q}$ may be compactified by one 
of the following
\begin{enumerate}
\item A cap trajectory $u\in \Nn^\circ(\xx,\zz,s;p)$
\item A triple $(u^-,v,u^+)$ with
$u^+ \in \Nn_{\rho^+}(\xx,\yy)$, $v \in \Nn^\U(\yy,\yyy;p)$ and $u^- \in 
	\Nn_{\rho^-}(\yyy,\zz)$
\item An unparametized gradient trajectory 
$u \in \widehat{\Nn}(\xx,\zz,s)$ and element $(R,v)$ with $v \in 
	\Nn_{\mu_R}(\xx,\zz;p)$
\item An element $(R,v)$ with $v \in \Nn_{\mu_R}(\xx,\zz;p)$
and an unparametized gradient trajectory $u \in \widehat{\Nn}(\xx,\zz,s)$
\end{enumerate}

We must show that in each of the above cases the relevant trajectories are 
thin.  In fact this follows from the conjugation by $p^+$ and $p^-$ and 
that $||H||< \sigma/2.$ Moreover that the limit trajectories are thin in 
cases (3) and (4) follows from the proof of Proposition 8 in 
\cite{chekanov:li98}. To economize the exposition we focus on cases (1) 
and (2).

Let us begin with (1).  The asymptotic limits assure us 
that $\xx = g_x \cdot \xxx$ and that $\zz = g_z \cdot \zzz$ with
$$
-\int (g_z \# \tilde{z'})^*\omega \leq 0 \ \textup{and} \ \int (g_x \# 
	\tilde{x'})^* \omega \leq 0.
$$
Hence the length of $u$ satisfies
\begin{eqnarray*}
\length(u) &=& \A_{sH}(\zz) - \A_{sH}(\xx) \\
&=& -\int \tilde{z}^*\omega - s\int H(z,t) \ dt + \int \tilde{x}^*\omega 
	+ s\int H(x,t) \ dt\\
&=& -\int (g_z\#\tilde{z'})^*\omega - s\int H(z,t) \ dt + \int 
	(g_x\#\tilde{x'})^*\omega + s\int H(x,t) \ dt\\
&\leq& s (\int H(x,t) - H(z,t) \ dt) \\
&\leq& s (\int \max_{x}H(x,t) - \min_{x}H(x,t) \ dt) \\
&=& s ||H|| \\
&<& \delta
\end{eqnarray*}
which makes it thin.

Figure \ref{fig:schematic} illustrates the situation of case (2).

\setlength{\unitlength}{1.65cm}
\begin{figure}
\begin{center}
\begin{picture}(6,4)(-3,-2)
\put(-2,1.5){\circle*{0.06}}
\put(-2,-1){\circle*{0.06}}
\put(1.5,1.5){\circle*{0.06}}
\put(1.5,-1){\circle*{0.06}}
\qbezier(-2,1.5)(-0.25,2)(1.5,1.5)
\qbezier(-2,1.5)(-0.25,1)(1.5,1.5)
\qbezier(-2,-1)(-0.25,-1.5)(1.5,-1)
\qbezier(-2,-1)(-0.25,-0.5)(1.5,-1)
\qbezier(1.5,1.5)(1,0.25)(1.5,-1)
\qbezier(1.5,1.5)(2,0.25)(1.5,-1)
\put(-.25,1.45){$u^+$}
\put(1.45,0.25){$v$}
\put(-.25,-1.05){$u^-$}
\put(-2.6,1.5){$\xx$}
\put(-2.6,-1){$\zz$}
\put(1.6,1.5){$\yy$}
\put(1.6,-1){$\yyy$}

\multiput(-3,-1.75)(6,0){2}{\line(0,1){4.}}
\multiput(-3,-1.75)(0,4.){2}{\line(1,0){6}}
\end{picture}
\caption{$u^+,v,u^-$ and their respective asymptotic 
limits}\label{fig:schematic}
\end{center}
\end{figure}
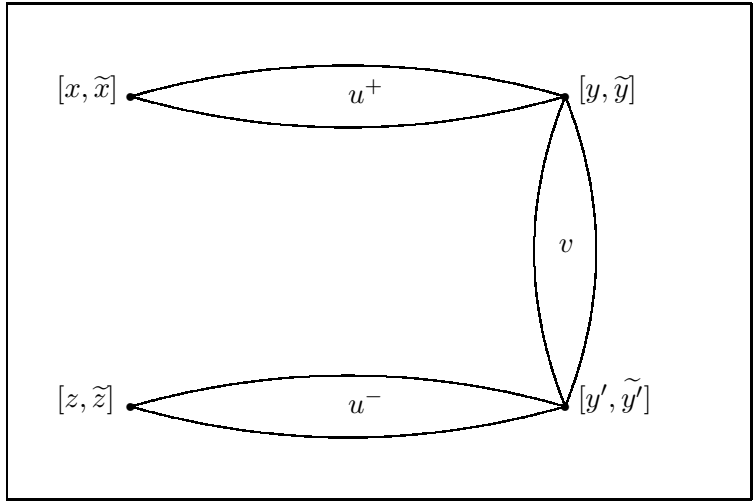

To complete the identification we must prove that each $u^+$ and $u^-$ is 
thin.  In other words we must prove $\length(u^+) \leq (1-s)\cdot b_+(H) + 
\delta$ and $\length(u^-) \leq (1-s)\cdot b_-(H) + \delta$.  
(Remember thin continuation trajectories satisfy Equation 
(\ref{eq:thin_continuation_def}).) Toward that 
end, consider the quantity $\length(u^+) + \length(u^-).$ By definition of 
length
\begin{eqnarray}
\length(u^+) + \length(u^-) &=& \A_{sH}(\xx) - \A_H(\yy) + \A_H(\yyy) - 
\A_{sH}(\zz) \nonumber \\
&=& \left( \int g_x^*\omega - s\int H(x,t) \ dt + \left( \A_H(\yyy) - 
	\A_H(\yy) \right) \right. \nonumber \\
&&\left. -\int g_z^*\omega + s\int H(z,t) \ dt\right). \label{eq:neg_terms}
\end{eqnarray}
The asymptotic conditions imply that 
$$
\int (g_x)^* \omega \leq 0 \ \textup{and} \ -\int (g_z)^*\omega \leq 0
$$
therefore (\ref{eq:neg_terms}) is smaller than
\begin{eqnarray}
&&s \left( \int H(z,t) \ dt) - \int H(x,t) \ dt \right) + \left( \A_H(\yyy) - 
	\A_H(\yy) \right).
	\nonumber
\end{eqnarray}
Since $v$ is a local cap trajectory from $\yy$ to $\yyy$ we know that 
$$
E(v) = \left( \A_H(\yyy) - \A_H(\yy) \right)
$$
and furthermore, because $v$ is local we have $E(v) \leq ||H|| < 
\delta/2.$  Therefore 
$$
\length(u^+) + \length(u^-) \leq \delta.
$$
Finally it follows that $u^+$ is thin:  Recall from Lemma 
\ref{lemma:cont_thin} that $\length(u^-) \geq -(1-s) \cdot b_+(H)$ (resp. 
$\length(u^+) \geq -(1-s)\cdot b_-(H)$) and so if $u^+$ were not thin then
\begin{eqnarray*}
\delta &=& -(1-s)\cdot b_+(H) + \delta + (1-s)\cdot b_+(H)\\
&<& \length(u^+) + \length(u^-) \leq \delta
\end{eqnarray*}
a contradiction.  Repeating the argument but interchanging $u^+$ and $u^-$ 
proves that $u^-$ is also thin.  In short the terms in Equation 
(\ref{eqn:to_verify}) are in one-one correspondence with the ends of the 
compact one-dimensional manifold $\mathcal{Q}$.  This implies $h^p$ is 
indeed a $\lambda-$homotopy.

\section{The Proof of Theorem \ref{thm:B} and Corollary \ref{corAB:corA}}

Theorem \ref{thm:B} is a consequence of Theorem \ref{thm:A} and
Proposition \ref{prop:A}.  Let $\mo$ be a closed symplectic manifold of 
dimension $2n$ and consider the product symplectic manifold
$(M\times M,-\omega \oplus \omega)$. The diagonal, $\Delta$, is a 
Lagrangian submanifold of dimension $2n = dim(M).$  In this case 
$\sigma(\Delta; M\times M, -\omega \oplus \omega) = \sigma(M;\omega)$.  
Fix a Darboux-Weinstein chart $(\Phi,\U)$ containing $\Delta.$  
Let the almost complex structure on $M\times M$ be of the form $\overline 
J=-J_{2t-1}\oplus J_{1-2t}$ 

Throughout this section consider the Hamiltonian deformation of the diagonal 
corresponding to the isotopy $t \mapsto \textup{graph}\ t \ d\Sp,$ and 
assume 
that the Hamiltonian function $H = H^\phi$ satisfies $||H|| = \osc(\Sp) < 
\sigma(M, \omega)/2$.  
Recall that the Hamiltonian function $\Hp$ is defined by the equation
$$
\Phi ( \Delta_{\phi_{\Hp}^t} ) = \ \textup{graph}\ t\ d\Sp
$$
and that the Hamiltonian vector field of $H$ has 
only fixed periodic orbits
$$
\textup{Per}(X_H) \longleftrightarrow \Delta \cap \Dp.
$$
For any non-degenerate Hamiltonian diffeomorphism we have 
$A(\phi;1) \leq \gamma (\phi) \leq ||\phi||_{Hofer} \leq ||H|| = \osc (\Sp).$  
Examination of Oh's proof in the $C^1-$small case makes clear that the 
crucial step of the proof is to produce a Floer trajectory connecting the 
maximum and minimum of the quasi-autonomous Hamiltonian $H$.  The trajectory we 
find is in fact local.  A curve in $M$ is called local if it is in the 
image of the association of Lemma \ref{le:open_closed}  Once this is 
established we go on to show $\osc(\Sp) \leq A(\phi;1)$, completing the 
proof.

Let $F = \Fp:(M \times M) \times [0,1] \to \R$ be the Hamiltonian function 
corresponding to $\Hp$, $\Fp(p,q,t) = \Hp(q,t).$  Set $y^\pm = 
(x^\pm,x^\pm) \in M \times M.$ 
\begin{lemma}\label{le:open_closed} There is an energy preserving one-one 
correspondance between the moduli spaces $\Nn^\U(\Delta, \Fp;\overline 
J; (q,q))$ and $\M^{local}(\Hp;J;q)$. 
\end{lemma}
\proof By definition any curve $v \in \Nn^\U(\Delta, 
\Fp;\overline{J};(q,q))$ with 
$$
\lim_{\tau \to \pm \infty} v(\tau,\cdot) = [y^\pm,\widetilde{y^\pm}] \ 
\textup{and} \ \widetilde{y^-} \#v \simeq \widetilde{y^+}
$$
satisfies
$$
\left\{
\begin{array}{l}
\begin{displaystyle}\dvdtau + \overline{J}(v,t)\left(\dvdt - X_F(v)\right) 
	= 0\end{displaystyle} \\ 
v:\R \times [0,1] \to M\times M \\
v(\tau,0),v(\tau,1) \subseteq \Delta \\
v(0,0) = (q,q) \in \Delta\\
\textup{image}(v) \subset \U 
\end{array}
\right.
$$
Integration by parts yields $E_{\overline J}(v) = \osc(\Sp).$  Indeed 
\begin{eqnarray*}
E_{\overline{J}}(v) &=& \int_{-\infty}^{+\infty} \left| \dvdtau \right|^2 
\ d\tau \\
&=& \int v^*(-\omega \oplus \omega) + ||\Fp||\\
&=& \int v^*(-\omega \oplus \omega) + ||\Hp||.
\end{eqnarray*}
Since the image of $v$ is contained in $\U$ the first term can be written
\begin{eqnarray*}
&&\int (\Phi \circ v)^*(d \Lambda_{can}) = 
\int_{\partial (\R \times [0,1])} (\Phi \circ v)^*\Lambda_{can}.
\end{eqnarray*}
By the boundary conditions this integral vanishes.  Hence 
$$
E_{\overline J}(v) = ||\Hp|| = \osc(\Sp).
$$

We glue the trajectory $v=(v_1,v_2)$ into a curve to $M$, 
$u:\R \times S^1 \to M$, as follows
$$
u(\tau,t) = \left\{
		\begin{array}{ll}
		v_1(2\tau, 2t-1) & 0 \leq t \leq 1/2\\
		v_2(2\tau, 1-2t) & 1/2 \leq t \leq 1
		\end{array}
		\right.
$$
The map $u$ is smooth by elliptic regularity and satisfies
$$
\left\{
\begin{array}{l}
\begin{displaystyle}\dudtau + J(u,t)\left( \dudt - X_H(u) \right) = 
	0\end{displaystyle} \\
u(\tau,t) = u(\tau, t+1) \\
\begin{displaystyle}\lim_{\tau \to \pm \infty} u(\tau, \cdot) = 
	x^\pm(\cdot)\end{displaystyle} \\
u(0,0) = p. 
\end{array}
\right.
$$
A straightforward chain rule calculation shows that $E_J(u) = E_{\overline
J}(v)=\osc(\Sp).$ \qed

Now we come to the main step in the proof of Theorem B.  Compare with 
Proposition 9.6, \cite{oh:si05}.

\begin{prop}\label{prop:local_curve}There exists an element $u \in 
\M^{local}([x^-,w^-], [x^+, w^+];p).$  Moreover any 
local curve with these asymptotics satisfies $E_J(u) = \osc(\Sp)$.
\end{prop}
\proof
We remind the reader (see Lemma 7, \cite{chekanov:li98} and Theorem 4.7, 
\cite{oh:fc96}) that when $s$ is sufficiently
small the thin part of the Floer homology satisfies 
$$HF^\circ_{2n-j}(\Delta, sF; \Z_2) \cong H_j(\Delta; \Z_2).$$  
Furthermore the thin part of the Floer cap action 
$[\Delta]\cap_s^\circ(\cdot)$ agrees with the standard cap action of 
$[\Delta]$ and this map is an isomorphism.  Let $j = 2n$ 
and $0$ and consider the following diagram
\begin{equation}\label{eq:commutative}
	\begin{array}{l}
\xymatrix{
  HF^\circ_{0}(\Delta, sF; \Z_2) \ar[r]^{\ \ \ \cong} 
\ar[d]_{[\Delta]\cap^\circ_s(\cdot)}^\cong \ar[d] &
  H_{2n}(\Delta; \Z_2)  \ar[d]_\cong^{[\Delta]\cap(\cdot)}\\	
  HF^\circ_{2n}(\Delta, sF; \Z_2) \ar[r]^{\ \ \ \cong} & H_0(\Delta;\Z_2)}
	\end{array}
\end{equation}

The thin cap action by the diagonal, the map $[\Delta] 
\cap^\circ_s(\cdot)$ in Equation (\ref{eq:commutative}), shifts the 
grading by the dimension of $\Delta.$

Because $\Delta$ is orientable $H_0(\Delta; \Z_2) \cong \Z_2$ and
$H_{2n}(\Delta; \Z_2) \cong \Z_2$ and so the 
thin part of the cap product on the Floer
homology $HF(\Delta, sF)$ is actually an
isomorphism on $\Z_2$ $$[\Delta]\cap^\circ_s(\cdot) :
\Z_2 \stackrel{\cong}{\to} \Z_2.$$ As such this map must be the
identity.  Hence it
commutes with the projection $p^+:CF(\Delta,sF) \to 
CF^+(\Delta;sF)\oplus CF^\circ(\Delta,sF)$ and 
Proposition \ref{prop:A} implies in
combination with Theorem \ref{thm:A} that the composition $\Phi^- \circ
[\Delta]\cap ^\U(\cdot) \circ \Phi^+$ is non-zero.  Because all three maps
are linear this means that the local Floer cap product
$[\Delta]\cap^\U(\cdot)$ is non-zero and therefore the moduli space 
$\Nn^\U
([(x^-,x^-),w^-], [(x^+,x^+), w^+];\Fp;\overline{J};p)$ must contain an 
odd 
number
of elements, where $x^\pm$ are respectively the maximum and minimum points
of $F$ and $w^\pm$ are any bounding disks. In particular it
is not empty.  Let us denote this element $v$.  By Lemma 
\ref{le:open_closed}, there exists $u \in \M^{local}([x^-,w^-], 
[x^+,w^+];q)$ and $E_J(u) = \osc(\Sp).$ 

\qed

\noindent \textit{Proof of Theorem B}

We now know that any local curve $u \in \M^{local}(H^\phi;p)$ satisfying 
$$
\lim_{\tau \to \pm \infty}u(\tau,\cdot) = x^\pm
$$
has energy equal to $\osc(\Sp)$ and this collection is not empty.  

On the other hand if $u' \in 
\M([x^-,\widetilde{\xm}],[x^+,\widetilde{\xp}];H^\phi;p)$ is any 
non-local curve 
then we have
\begin{eqnarray*}
E_J(u') &=& E_{J'}(v')
\end{eqnarray*}
where $v'(\tau,t) = \left(\phi^t_H\right)^{-1} \left( u'(\tau,t) \right)$.
One may check that the curve $v'(\tau,t) = \left(\phi^t_H\right)^{-1} 
\left( u'(\tau,t) \right)$ satisfies the equation
$$
\dvdtau' + J'\dvdt' = 0
$$
where $J' = (\phi^t)^*J$.  Since the periodic orbits of $X_{\Hp}$ are 
all constant, the glued curve $v'\#v$, where 
$v = v(\tau,t) = \left(\phi^t_H\right)^{-1} \left( u(\tau,t) \right)$ and 
$u$ is local curve we constructed in Proposition 
\ref{prop:local_curve} produces a 
non-constant pseudo-holomorphic disk with 
boundary on the diagonal $\Delta$.

Therefore 
$$
\sigma(M,\omega) < E(v'\#v) = E(v') + ||\Hp||.
$$  
Since we assume that $$\osc(\Sp) = ||\Hp|| < \sigma(M,\omega)/2,$$
it must be that $E_{J'}(v') > \sigma/2.$
In short, $E_J(u') \geq ||\Hp||.$

Therefore $A(\phi,J_0;J;[\xm,\widetilde{\xm}], [\xp,\widetilde{\xp}];p) = 
||\Hp|| = \osc(\Sp).$  
Since $\Hp$ is quasi-autononous we then find $A(\phi,J_0;J;1;p) = 
\osc(\Sp).$  
Therefore 
the supremum over $p$, $J$ and $J_0$ is greater than the oscillation, in 
other words we have proven
$$A(\phi;1) \geq \osc(\Sp).$$ 
This proves Theorem \ref{thm:B}.\qed

Before proving Corollary \ref{corAB:corA}, we prove a simple estimate 
between the $C^1$ norm of a Hamiltonian flow to the Hofer 
length of a generating Hamiltonian function.  Let $J$ be any 
$\omega-$compatible almost complex structure and denote $g_{J}$ the 
corresponding metric.  With respect to the metric $g_{J}$, the
$C^1$ norm of the isotopy $\phi_H$ is
\begin{eqnarray}
||\phi_H||_{C^1(M\times [0,1])} &:=& \sup_{x,t} 
||\left(T\phi^t_H\right)_{|x}|| + \sup_{x,t} ||\dot{\phi^t_H}(x)||\\
&=& \sup_{x,t} 
||\left(T\phi^t_H\right)_{|x}|| + \sup_{x,t} ||X_H(x,t)||.\nonumber
\end{eqnarray}
\begin{lemma}\label{lemma:esti} For any Hamiltonian isotopy $t \mapsto 
\phi^t_H$ 
\begin{equation}
||H|| \leq \textup{diam}_{g_J}(M) 
\cdot 
||\phi_H||_{C^1(M\times [0,1])}.
\end{equation}
\end{lemma}
\proof  Let $d$ denote the Riemannian distance on $M$ and $x^+_t$ and 
$x^-_t$ be the points in $M$ where extreme values of $H(\cdot,t)$ occur.
$$\max_{x\in M} H(x,t) = H(x^+_t,t) \ \textup{and} \ \min_{x\in M}H(x,t) = 
H(x^-_t,t)$$  
For $\epsilon > 0$ choose a path $\sigma_t =\sigma_t(s) : [0,1] \to M$ 
with 
$$
\sigma_t(0) = x^-_t, \ \sigma(1) = x^+_t \ \textup{and} \ \int_0^1 
|\dot {\sigma_t}|_{g_J} \ ds < \epsilon + d(x^-_t,x^+_t).$$
Next
\begin{eqnarray}
H(x^+_t,t) - H(x^-_t,t) &=& \int_0^1 \frac{d}{ds} H(\sigma(s),t) \ ds 
\nonumber\\
&=& \int_0^1 g_J(\nabla H(\sigma, t), \dot \sigma) \ ds. 
\label{eq:est_grad}
\end{eqnarray}
By the Cauchy-Schwarz inequality (\ref{eq:est_grad}) is smaller than
$$||\nabla H||_{L^\infty} \cdot ||\dot 
\sigma||_{L^1} \leq ||\nabla H||_{L^\infty} (d(x,y) + \epsilon).
$$
Since $\epsilon>0$ is arbitrary it follows that 
$$||H||_{Hofer} \leq ||\nabla H||_{L^\infty}\cdot d(x^-_t,x^+_t)
\leq ||\nabla H||_{L^\infty} \cdot \textup{diam}_{g_J}(M).$$
Finally, the compatible almost complex strucure $J$ is a $g_J$ isometry, 
hence $||\nabla H|| = ||X_H||$.  This proves the Lemma. \qed\\

\noindent \textit{Proof of Corollary \ref{corAB:corA}.}
A priori we know that
\begin{equation}\label{eq:apriori}
\gamma(\phi) \leq ||\phi||_{Hofer} \leq ||H^\phi||.
\end{equation}

Let $\epsilon > 0$ be such that $\osc(\Sp) + \epsilon < \sigma\mo /2$ 
and $\psi^t$ be a non-degenerate Hamiltonian isotopy with $$||\phi^t 
\circ (\psi^t)^{-1}||_{C^1(M\times [0,1])} < 
\frac{\epsilon}{\textup{diam}_{g_J}(M)}.$$
Theorem \ref{thm:B} implies $\gamma(\psi) = ||\psi||_{Hofer} = ||H^\psi||$ 
and Lemma \ref{lemma:esti} 
implies $$||H^\phi \# \overline{H^\psi}|| < \epsilon.$$
By the triangle inequality,  
\begin{eqnarray*}
||H^\phi|| &\leq& ||H^\phi \# \overline{H^\psi}|| + ||H^\psi||\\
&<& \epsilon + ||H^\psi||\\
&=& \epsilon + \gamma(\psi) \\
&\leq& \epsilon + \gamma(\psi \circ \phi^{-1}) + \gamma(\phi).
\end{eqnarray*}
By the symmetry property of the spectral norm we 
next see that the above equals
\begin{eqnarray*}
&=& \epsilon + \gamma(\phi \circ \psi^{-1}) + \gamma(\phi)\\
&\leq& \epsilon + ||H^\phi \# \overline{H^\psi}|| + \gamma(\phi).
\end{eqnarray*}
Hence $$0\leq ||H^\phi|| - \gamma(\phi) < 2\epsilon.$$
Letting $\epsilon$ approach zero proves $||H^\phi|| = \gamma(\phi).$  
Applying Equation (\ref{eq:apriori}) a second time proves 
$$\gamma(\phi) = ||\phi||_{Hofer} = ||H^\phi||.$$ \qed

\smallskip
\begin{center}
\textsc{Concluding Remarks}
\end{center}
\smallskip

It can be rather difficult to compute the constant $\sigma \mo$.  For the
two-sphere with total surface area equal to $4\pi$ then in fact
$\sigma(S^2 \times S^2, -\omega \oplus \omega)/2 = 2\pi.$ The Hamiltonian
function which generates the half-turn rotation about a diameter has Hofer
length $2\pi$.  In this sense the constant $\sigma/2$ seems to be precise.  
This is plausible as rotation by more than half the way around a diameter
is not Hofer (nor spectral) length minimizing:  one can generate the same
diffeomorphism by rotation in the opposite direction.  However it is not
yet know if this situation is engulfable.  In some cases the Darboux
neighborhood of the diagonal can be quite large.  For the sphere the 
author speculates the largest possible Darboux neighborhood to be the 
compliment of the anti-diagonal in $S^2 \times S^2$. 

\bibliographystyle{amsplain}
\bibliography{length_min.bib}
\end{document}